\documentclass[11pt]{article}
\usepackage{amsmath}
\usepackage{latexsym}

\newtheorem{theorem}{Theorem}[section]

\newtheorem{conjecture}[theorem]{Conjecture}

\newtheorem{definition}[theorem]{Definition}

\DeclareMathOperator{\Pyr}{Pyr}
\DeclareMathOperator{\Prism}{Prism}
\DeclareMathOperator{\Bipyr}{Bipyr}

\begin{document}

\title{The $cd$-Index: A Survey}

\author{Margaret M. Bayer \\
Department of Mathematics\\ University of Kansas\\
       Lawrence, KS, USA \\ bayer@ku.edu}

\maketitle

\begin{abstract}
This is a survey of the $cd$-index of Eulerian partially ordered sets.  
The $cd$-index is an encoding of the numbers of chains, specified by
ranks, in the poset.  It is the most efficient such encoding, 
incorporating all the affine relations on the flag numbers of Eulerian
posets.  Eulerian posets include the face posets of regular CW spheres
(in particular, of convex polytopes), intervals in the Bruhat
order on Coxeter groups, and the lattices of regions of oriented matroids.
The paper discusses inequalities on the $cd$-index, connections with
other combinatorial parameters, computation, and algebraic approaches.
\end{abstract}

\section{Early History}
The history of the $cd$-index starts with the combinatorial study of convex
polytopes.
Over one hundred years ago Steinitz proved the characterization of the
face vectors of 3-dimensional polytopes \cite{Steinitz}.
Interest in the number of faces of convex polytopes in higher dimensions
grew with the development of linear programming from the 1950s on.  
While the affine span of the face vectors of $d$-dimensional polytopes
is known to be given just by the Euler equation, a full
characterization of the
face vectors of polytopes of dimensions 4 and higher still eludes us. 
The major breakthrough on this question came through the definition of
the Stanley--Reisner ring and the discovery of the connection between
polytopes and toric varieties.  
As a result the face vectors of simplicial polytopes were characterized
by Billera and Lee~\cite{Bil-Lee} and Stanley~\cite{Stanley-simplicial}, 
following a conjecture of McMullen~\cite{McMullen}.

In the 1970s Stanley, in studying balanced complexes and ranked partially
ordered sets \cite{Stanley-balanced}, broadened
the focus from counting the elements of each rank to counting chains in
the poset with elements from a specified rank set.  
Bayer and Billera applied this perspective to convex polytopes, and
initiated the study of ``flag vectors'' of polytopes \cite{Bayer-Bil}.
They found the ``generalized Dehn--Sommerville equations,''
a complete set of equations defining the affine span of the
flag vectors of $d$-polytopes; the dimension turns out to be a Fibonacci 
number.  As with face vectors, a complete characterization of flag 
vectors of polytopes is unknown.

Shortly after the proof of  the generalized Dehn--Sommerville equations,
Fine (see \cite{Bayer-Klapper}) found a way to encode the flag vectors in the 
most efficient way, with the $cd$-index.  
The generalized Dehn--Sommerville equations apply
not just to convex polytopes, but to all Eulerian posets.
Likewise, the $cd$-index is defined for Eulerian posets, which include
the face posets of regular CW spheres, intervals in the Bruhat order on finite
Coxeter groups, and the lattices of regions of oriented matroids.  

There are two main issues for research on $cd$-indices.  One is the
question of the nonnegativity of the coefficients, or, more generally,
inequalities on the $cd$-index, for Eulerian posets or for particular
subclasses.
The other (related) issue is the combinatorial interpretation of the 
coefficients, either directly in terms of the poset or in terms of
other combinatorial objects.
In the last thirty-five years, much research has been carried out on
these issues.

\section{Basic Definitions}\label{Basic}
\begin{definition}{\em 
An {\em Eulerian poset} is a graded partially ordered set such that each
(nonsingleton) interval $[x,y]$ in the poset has an equal number of elements 
of even
and odd rank.  }\end{definition}

\begin{definition}[from \cite{Stanley-Euler}]{\em 
A finite {\em regular CW complex} is a finite collection
of disjoint open cells $\sigma$ in Euclidean space such that each 
$\sigma$ is homeomorphic to an open ball of some dimension $n$ and
its boundary is homeomorphic to a sphere of dimension $n-1$, which
is the union of lower dimensional cells.
If the complex is homeomorphic to a sphere, it is called a 
{\em regular CW sphere}.}\end{definition}

(For a precise definition that does not assume Euclidean space, see
\cite{Hersh}.)

Bj\"{o}rner~\cite{Bjorner}  showed that the face posets of CW complexes are
exactly the posets with a unique minimum element $\hat{0}$  and at least one
other element, and for which the order complex of every open interval
$(\hat{0},x)$ is homeomorphic to a sphere.
(The {\em order complex} of a poset is the simplicial complex whose faces are 
the chains of the poset.)

The set of closed cells of a regular CW sphere, ordered by inclusion,
along with the empty set and an adjoined maximum element, forms an Eulerian
poset.  A convex polytope is a regular CW sphere. 
The face lattice of an $n$-dimensional polytope is a rank $n+1$ Eulerian 
poset.
(In what follows we generally do not distinguish between a polytope and its
face lattice.)

\begin{definition}{\em 
For $P$ a graded poset of rank $n+1$ and $S\subseteq[n]=\{1,2,\ldots, n\}$,
the {\em $S$-flag number} of $P$, denoted $f_S(P)$, is the number of chains 
$x_1\prec x_2 \prec \cdots \prec x_s$ of $P$ for which 
$\{\mbox{rank}(x_i): 1\le i\le s\} = S$.  (By convention, 
$f_\emptyset(P)=1$.)  The {\em flag vector} of $P$ 
is the length $2^n$ vector $(f_S(P))_{S\subseteq[n]}\in {\bf N}^{2^n}$.
}\end{definition}

When the poset is the face lattice of a polytope, the indexing set is typically
shifted to represent dimensions of the faces, rather than ranks in the poset.
The restriction of the flag vector to terms indexed by singleton sets is
the {\em $f$-vector}  (or {\em face vector}) of $P$.

Sommerville~\cite{Sommerville} proved the Dehn--Sommerville equations for
$f$-vectors of simplicial polytopes by applying Euler's formula to
each interval in the face lattice.  The same method gives equations on
flag vectors, known as the {\em generalized Dehn--Sommerville equations}.

\begin{theorem}[\cite{Bayer-Bil}]
The affine dimension of the flag vectors of rank $n+1$ Eulerian posets is
$F_n-1$, where $(F_n)$ is the Fibonacci sequence (with $F_0 = F_1 = 1$).
The affine hull of the flag vectors is given by the equations
\begin{equation*}\sum_{j=i+1}^{k-1} (-1)^{j-i-1}f_{S\cup\{j\}}(P) = 
(1-(-1)^{k-i-1})f_S(P),\end{equation*}
where $i\le k-2$, $i,k\in S\cup\{0,n+1\}$, and 
$S\cap \{i+1,\ldots, k-1\}=\emptyset$.
\end{theorem}

This affine space is spanned by the flag vectors of convex polytopes.
Bases of polytopes were given by Bayer and Billera~\cite{Bayer-Bil}
and by Kalai~\cite{Kalai}.

In the mid-1980s, Jonathan Fine discovered a compact way to represent these
equations.  
To see this, we first need the transformation to the flag $h$-vector.

\begin{definition}{\em
\begin{sloppypar}
Let $P$ be a rank $n+1$ Eulerian poset with flag vector 
$(f_S(P))_{S\subseteq[n]}$.  The {\em flag $h$-vector} of $P$ is the vector
$(h_S(P))_{S\subseteq[n]}\in{\bf N}^{2^n}$, where
\begin{equation*}h_S(P) = \sum_{T\subseteq S} (-1)^{|S\setminus T|} f_T(P).
\end{equation*}
\end{sloppypar}
}\end{definition}

This transformation is invertible: $f_S(P) = \sum_{T\subseteq S}  h_T(P)$.
The flag $h$-vector has algebraic meaning through the Stanley--Reisner ring
of the order complex.
For a convex polytope (and more generally, a balanced Cohen--Macaulay
complex), the entries in the flag $h$-vector are nonnegative
\cite{Stanley-balanced}.

The flag $h$-vector can be represented by  a polynomial in noncommuting 
variables $a$ and $b$.  Associate with  $S\subseteq[n]$ the monomial
$u_S = u_1u_2\cdots u_n$, where $u_i = a$ if $i\not\in S$ and 
$u_i = b$ if $i\in S$.  Then the $ab$-polynomial is 
$\Psi_P(a,b) = \sum_{S\subseteq[n]} h_S u_S$.
Here is an equivalent formulation for the $ab$-polynomial: associate to
each chain $x_1 \prec x_2 \prec \cdots \prec x_s$ of $P$ with rank set 
$S$ the monomial
$w_1w_2\cdots w_n$, where $w_i = a-b$ if $i \not\in S$ and $w_i=b$
if $i\in S$.  Then $\Psi_P(a,b)$ is the sum of these monomials over
all chains of $P$.

Fine's inspiration was to see that when $P$ is a convex polytope, the
$ab$-polynomial can be written as a polynomial with integer coefficients
in the noncommuting variables $c$ and $d$, where $c=a+b$ and $d=ab+ba$.

\begin{definition}{\em
Let $P$ be a rank $n+1$ poset.  The {\em $cd$-index} of $P$ is the
polynomial $\Phi_P(c,d)$ such that $\Phi_P(a+b,ab+ba)= \Psi_P(a,b)$,
if such a polynomial exists.
}\end{definition}

The $cd$-index of a rank $n+1$ poset is considered a homogeneous 
polynomial (in noncommuting variables)
of degree $n$ by assigning degree 1 to $c$ and degree 2 to $d$.
A straightforward recursion shows that the number of $cd$-words
of total degree $n$ is the Fibonacci number $F_n$.
It is easy to see from the definition that the $cd$-index of the
dual of a poset (reverse the order relation) is obtained from the
$cd$-index of the poset by reversing all the $cd$-words.

\begin{theorem}[\cite{Bayer-Klapper}]
Let $P$ be a graded poset.  Then $P$ has a $cd$-index 
if and only if the flag $f$-vector of $P$ satisfies
the generalized Dehn--Sommerville equations.
In this case the coefficients of the $cd$-index are integers.
\end{theorem}

We will sometimes refer to the affine space of coefficients of the
$cd$-words of fixed degree as the {\em generalized Dehn--Sommerville
space}.

The definition of the $cd$-index gives a way of computing it from
the flag \mbox{$h$-vector}, and hence from the flag $f$-vector.  
Here are formulas
for several low ranks.  (Meisinger's dissertation~\cite{Meisinger} has
many useful tables, including flag number formulas for the $cd$-index 
up through rank~9.)
\begin{description}
\item[rank 3] $cc+(f_1-2)d$
\item[rank 4] $ccc+(f_3 - 2)cd+(f_1-2)dc$
\item[rank 5] $cccc+(f_4-2)ccd+(f_2-f_1)cdc+(f_1-2)dcc+(f_{13}-2f_3-2f_1+4)dd$
\end{description}

Table 1 gives the $cd$-indices of some familiar polytopes.

\begin{table}[h]
\caption{Some $cd$-indices}
\begin{center}
\begin{tabular}{|l|l|}
\hline
$n$-gon & $cc + (n-2)d$ \\
\hline
tetrahedron & $ccc+2cd+2dc$ \\
\hline
cube & $ccc+4cd+6dc$ \\
\hline
octahedron & $ccc+6cd+4dc$ \\ 
\hline
4-simplex & $cccc+3ccd+5cdc+3dcc+4dd$ \\
\hline
4-cube & $cccc+6ccd+16cdc+14dcc+ 20dd$ \\
\hline
\end{tabular}
\end{center}
\end{table}

Fine believed that the coefficients of the $cd$-index of a convex
polytope are always nonnegative.  This appears as a more general
conjecture in \cite{Bayer-Klapper}.

\begin{conjecture}
The coefficients in the $cd$-index of every regular CW sphere are
nonnegative.
\end{conjecture}
The conjecture turns out to be true.
In the next section we will look at the great body of work addressing
this conjecture.

Note that while Fine did not publish anything about the $cd$-index, 
his calculations involving the $cd$-index inspired his work on an
alternative approach to flag vectors \cite{Fine,Fine-arxiv}.

Stanley \cite{Stanley-flag}
noted that it is sometimes useful to write the $cd$-index
as a polynomial in $c$ and $e^2$, where
$e=a-b$ and thus $e^2 = c^2 - 2d$.
Purtill~\cite{Purtill} showed that
if $P$ is a convex polytope, then $\Phi_P(c,d)$ can be written as a polynomial
in the noncommuting variables $c$, $d$, and $-e^2 = 2d-c^2$ with 
nonnegative coefficients. 

Next, let us consider another important parameter for convex polytopes.
The {\em $h$-vector} of a simplicial polytope is the result of a certain
linear transformation on the $f$-vector.  This transformation
was noted by Sommerville~\cite{Sommerville}, but its significance was
not understood for decades.  For simplicial polytopes, the $h$-vector
has interpretations in shellings of polytopes, the Stanley--Reisner
ring and the toric variety associated with the polytope.
Unfortunately, the particular transformation from $f$-vector to $h$-vector
does not give a meaningful vector if the polytope is not simplicial.
This problem was resolved by consideration of the toric variety associated
with a nonsimplicial, rational polytope, and, in~particular, its 
intersection homology.  Stanley~\cite{Stanley-h} gave the definition
as follows.  In the case of simplicial polytopes, this $h$-vector
specializes to the aforementioned $h$-vector.  In the general case it 
is sometimes referred to as the ``toric $h$-vector.''  

The recursive definition uses the following notational conventions.
For an Eulerian poset $P$, write $\hat{0}$ for the unique minimal
element and $\hat{1}$ for the unique maximal element.  Denote by $P^-$
the poset $P\setminus \{\hat{1}\}$.
We use interval notation in a poset; in particular
$[\hat{0},t) = \{s\in P: \hat{0}\preceq s\prec t\}$.  The rank of
an element $t$ of $P$ is denoted~$\rho(t)$.  Finally, $k_{-1} = 0$.

\begin{definition}\label{torich}{\em
Families of polynomials $f$ and $g$ in a single variable $x$ are defined by 
the following rules:
\begin{itemize}
\item $f(\emptyset,x) = g(\emptyset,x) = 1$
\item If $P$ is an Eulerian poset of rank $n+1\ge 1$, and if
      $f(P^-,x) = \sum_{i=0}^n k_ix^i$, then 
      $g(P^-,x) = \sum_{i=0}^{\lfloor n/2\rfloor} (k_i - k_{i-1})x^i$.
\item If $P$ is an Eulerian poset of rank $n+1\ge 1$, then

      $f(P^-,x) = \sum_{t\in P^-} g([\hat{0},t),x)(x-1)^{n-\rho(t)}$
\end{itemize}
When $\displaystyle f(P^-,x)= \sum_{i=0}^n k_ix^i$, let $h_i = k_{n-i}$, and 
call $(h_0,h_1,\ldots, h_n)$ the (toric) $h$-vector of  $P$.
}\end{definition}

In the case where $P$ is the face lattice of a rational convex
$n$-polytope, this $h$-vector has many of the same properties as the 
$h$-vector of a simplicial polytope; in particular, it is nonnegative
($h_i\ge 0$) and symmetric ($h_i=h_{n-i}$).

\begin{theorem}
For Eulerian posets, the toric $h$-vector is the image of a linear 
transformation of the flag vector.  Equivalently, it is the image of a
linear transformation of the $cd$-index.
\end{theorem}

A recursive proof of this proposition was given in \cite{Bayer-Klapper}.
Explicit formulas for the $cd$-index--$h$-vector map were given in 
\cite{Bayer-Ehr}.
Hetyei~\cite{Hetyei-short} was able to simplify the $cd$-index--$h$-vector 
relation by introducing the ``short toric polynomial.''
The paper~\cite{Bayer-Ehr}  also gave a combinatorial approach involving 
lattice paths, due to Fine.
A method for computing the toric $h$-vector (as well as the $cd$-index)
of a convex polytope by sweeping a hyperplane through the polytope was given by
Lee in \cite{Lee-sweeping}.

\section{Inequalities}
\subsection{Nonnegativity}\label{Nonneg}

As mentioned before, Fine believed that the $cd$-index has nonnegative
coefficients for all polytopes, and Bayer and Klapper conjectured 
nonnegativity for all regular CW spheres.  Here is a sequence of 
results and conjectures on nonnegativity.  We say the $cd$-index of 
a poset is {\em nonnegative} when all its coefficients are 
nonnegative, and we denote this by $\Phi_P(c,d)\ge 0$.

\begin{theorem}[\cite{Purtill} Purtill 1993]\label{Purtill}
The $cd$-indices of the following polytopes are nonnegative.
\begin{itemize}
\item polytopes of dimension at most 5
\item simple polytopes
\item simplicial polytopes
\item quasisimplicial polytopes (all facets are simplicial)
\item quasisimple polytopes (all vertex figures are simple)
\end{itemize}
\end{theorem}

\begin{theorem}[\cite{Stanley-flag} Stanley 1994]
The $cd$-indices of $S$-shellable CW spheres are nonnegative.
In particular the $cd$-index of every convex polytope is nonnegative.
\end{theorem}

(The class of $S$-shellable CW spheres includes all convex polytopes.)

An important class of posets is the class of Cohen--Macaulay posets.
These are the posets whose order complexes have Cohen--Macaulay Stanley--Reisner
rings.
A poset is {\em Gorenstein*} if and only if it is Cohen--Macaulay and Eulerian.

\begin{conjecture}[\cite{Stanley-flag,Stanley-Euler}]

\mbox{}

\begin{itemize}
\item The $cd$-index of every Gorenstein* poset is nonnnegative.
\item The $cd$-index of every Gorenstein* lattice is coefficientwise greater
      than or equal to the $cd$-index of the Boolean lattice (simplex).
\end{itemize}
\end{conjecture}
Stanley~\cite{Stanley-flag} showed that if 
a Gorenstein* poset is also simplicial (all of its intervals up to an element
$x\ne \hat{1}$ are Boolean lattices), then its $cd$-index is nonnegative.
Moreover, he showed that the inequalities of this theorem would
imply all linear inequalities satisfied by the flag $f$-vectors of all
Gorenstein* posets, and all those satisfied by the smaller class of
$S$-shellable CW spheres.

In the most general case, Bayer~\cite{Bayer-signs} determined which 
$cd$-coefficients are bounded for all Eulerian posets.
\begin{theorem}\label{signs}

\mbox{}

\begin{enumerate}\item\label{th_pt1}
 For the following $cd$-words $w$, the coefficient of $w$ as a function
of Eulerian posets has greatest lower bound $0$ and has no upper bound:
\begin{enumerate}
\item $c^idc^j$, with $\min\{i,j\}\le 1$,
\item $c^idcd\cdots cdc^j$ (at least two $d$'s alternating with $c$'s, $i$ and
      $j$ unrestricted).
\end{enumerate}
\item\label{th_pt2} The coefficient of $c^n$ in the $cd$-index of every
      Eulerian poset is $1$.
\item\label{th_pt3}
For all other $cd$-words $w$, the coefficient of $w$ as a function of
Eulerian posets has neither lower nor upper bound.
\end{enumerate}
\end{theorem}

In particular, there are Eulerian posets with some negative $cd$-coefficients.
Easy examples in odd rank are obtained as follows.
Take two copies of the ``smallest'' rank~$n$
Eulerian poset, namely, the poset with two elements at each rank, all pairs
of elements of different ranks comparable. 
Identify the top and bottom elements of the two copies.  
For rank $2k+1$, the $cd$-index is
$2c^{2k}-(c^2-2d)^k$ \cite{Ehr-Rea-non}.

Before Stanley's conjecture was proved by Karu for all Gorenstein* posets,
there was progress on special cases.
Reading~\cite{Reading-nonneg}
used the proof of the Charney--Davis Conjecture for dimension 3
(Davis and Okun~\cite{Davis-Okun}) and a convolution formula to prove the 
nonnegativity of the coefficients of certain $cd$-words for 
Gorenstein* posets.
Novik~\cite{Novik} proved the nonnegativity of certain $cd$-coefficients 
for odd-dimensional simplicial complexes that are Eulerian and Buchsbaum
(a weakening of Cohen--Macaulay), in particular for odd-dimensional
simplicial manifolds.
Hetyei~\cite{Hetyei-polysph} constructed a set of polyspherical 
CW complexes having nonnegative $cd$-indices.  (The face posets of these
complexes are Gorenstein* posets.)
Hsiao~\cite{Hsiao} constructed an analogue of distributive lattices
having nonnegative $cd$-indices.  (These are Gorenstein* posets.)

Karu  pursued a proof of nonnegativity of the $cd$-index for complete
fans using methods from algebraic geometry, and was able to extend his proof
to the general Gorenstein* case.

\begin{theorem}[\cite{Karu-cd} Karu 2006]
The $cd$-index of every Gorenstein* poset is nonnegative.
\end{theorem}

A complete characterization of the $cd$-indices of Gorenstein* posets is
presumably beyond reach, but Murai and Nevo~\cite{Murai-Nevo-5} obtained
the characterization for rank~5.

\subsection{Monotonicity}\label{Monotonicity}
In this section we consider comparisons among the $cd$-indices of different
posets.

Billera, Ehrenborg and Readdy studied the $cd$-indices of zonotopes (polytopes
arising as the Minkowski sum of segments) and, more generally, of the
lattice of regions of oriented matroids.
In analogy to our notation for nonnegativity of the $cd$-index, we write
$\Phi_Q(c,d)\ge \Phi_P(c,d)$ to mean that the coefficient of each $cd$-monomial
in the $cd$-index of $Q$ is greater than or equal to the corresponding
coefficient in the $cd$-index of $P$.
\begin{theorem}[\cite{c2d}]
Let $Q_n$ be the $n$-dimensional cube, and $Q_n^*$ its dual poset.
\begin{itemize}
\item
If the rank $n+1$ poset $P$ is the lattice of regions of an oriented matroid,
then $\Phi_P(c,d)\ge \Phi_{Q_n^*}(c,d)$.
\item If $Z$ is an $n$-dimensional zonotope, then 
$\Phi_Z(c,d)\ge \Phi_{Q_n}(c,d)$.
\end{itemize}
\end{theorem}
In this context the $cd$-index can be modified to the $c$-$2d$-index, 
because the
coefficient of every word containing $k$ $d$s is a multiple of $2^k$.

Nyman and Swartz fixed the dimension and number of zones and found the 
zonotopes with minimum and maximum $cd$-indices.

\begin{theorem} [\cite{Nyman}]
For fixed $r$ and $n$, let ${\mathcal H}_L$ be an essential hyperplane arrangement
with underlying geometric lattice the rank $r$ near pencil with $n$ atoms, 
and let ${\mathcal H}_U$  be an essential hyperplane arrangement
with underlying geometric lattice the rank $r$ truncated Boolean lattice
with $n$ atoms. Let $Z_L$ and $Z_U$ be the zonotopes dual to 
${\mathcal H}_L$  and ${\mathcal H}_U$.
Then for any $r$-dimensional zonotope $Z$ with $n$ zones, 
\begin{equation*}\Phi_{Z_L}(c,d)\le \Phi_Z(c,d)\le \Phi_{Z_U}(c,d).\end{equation*}
\end{theorem}

Ehrenborg~\cite{Ehr-zono} gave additional inequalities for the $cd$-index
of zonotopes.

The following result of Billera and Ehrenborg is analogous to an inequality
on the toric $g$-vector of rational polytopes, conjectured by Kalai~\cite{Kalai}
and proved by Braden and MacPherson~\cite{Braden}.

\begin{theorem}[\cite{Bil-Ehr}]
Let $P$ be a polytope and $F$ a face of $P$.  Let $P/F$ be the polytope whose
face lattice is the interval $[F,P]$ of the face lattice of $P$.  Denote
the pyramid over a polytope $Q$ by $\Pyr(Q)$.
Then
\begin{itemize}
\item $\Phi_P(c,d)\ge \Phi_F(c,d)\cdot \Phi_{\Pyr(P/F)}(c,d)$
\item $\Phi_P(c,d)\ge \Phi_{\Pyr(F)}(c,d)\cdot \Phi_{P/F}(c,d)$
\end{itemize}
\end{theorem}

In particular, if $F$ is a facet of $P$, then $\Phi(P)\ge \Phi(\Pyr(F))$,
so among all polytopes having $F$ as a facet, the one with minimum $cd$-index
is the pyramid over~$F$.  Repeated application of this shows that the simplex
minimizes the $cd$-index among polytopes.

Billera and Ehrenborg were also able to show the upper bound theorem for
$cd$-indices of polytopes.
(For the upper bound theorem for $f$-vectors, see \cite{McMullen}.)
\begin{theorem}[\cite{Bil-Ehr}]
Let $P$ be an $r$-dimensional polytope with $n$ vertices, and let 
$C(n,r)$ be the cyclic $r$-polytope with $n$ vertices.  Then
\begin{equation*}\Phi_P(c,d)\le \Phi_{C(n,r)}(c,d).\end{equation*}
\end{theorem}

Ehrenborg and Karu proved a decomposition theorem for the $cd$-index of a
Gorenstein* poset, resulting in the following inequalities.

\begin{theorem}[\cite{Ehr-Karu}]
Let $B_n$ be the Boolean lattice of rank $n$.
\begin{itemize} 
\item If P is a Gorenstein* lattice of rank $n$, then 
      $\Phi_P(c,d)\ge \Phi_{B_n}(c,d)$.
\item If $P$ is a Gorenstein* poset, and $Q$ is a subdivision of $P$,
      then \linebreak $\Phi_Q(c,d)\ge \Phi_P(c,d)$.
\end{itemize}
\end{theorem}

In his master's thesis, Dornian proved the following.
\begin{theorem}[\cite{Dornian}]
Let $\Pi$ be the face poset of a simplicial $(d-1)$-sphere, let
$S_d(k)$ be a stacked $d$-polytope, and suppose each has a triangulation as a 
shellable $d$-ball with $k$ simplices.
Then ${\Phi_\Pi\le \Phi_{S_d(k)}}$.
\end{theorem}

\subsection{Other Inequalities}\label{other-ineq}
Stanley~\cite{Stanley-flag} showed that for each $cd$-word $w\ne c^n$
there is a sequence of Eulerian posets whose $cd$-indices (normalized)
tend to $w$.
This can be seen as a strengthening of the fact that coefficients of 
$cd$-words have no upper bound (Theorem~\ref{signs}).
Another proof of this by Bayer and Hetyei is in~\cite{BayHet-flag}, 
where some extreme rays
of the closed cone of flag $f$-vectors of Eulerian posets are given.

The nonnegativity of the $cd$-index can be translated into inequalities
on the flag $h$-vector and flag $f$-vector.  Some simpler inequalities
can also be extracted.  Stanley~\cite{Stanley-flag} considered the
comparison of two entries of the flag $h$-vector.  The result depends
on a function of sets that looks mysterious, but makes more sense when
visualizing how a $cd$-word (with $d$ of degree 2) ``covers'' an
interval of integers.  
For $S\subseteq [n]$, let 
$\omega(S) = \{i\in[n-1]: \mbox{exactly one of $i$ and $i+1$ is in $S$}\}$.
Stanley~\cite{Stanley-flag} showed that the following theorem follows from
the nonnegativity of the $cd$-index for all Gorenstein* posets.
\begin{theorem}
Let $S$ and $T$ be subsets of $[n]$.  The following are equivalent.
\begin{itemize}
\item $\omega(T)\subseteq\omega(S)$
\item For every Gorenstein* poset $P$ of rank $n+1$, $h_T(P)\le h_S(P)$.
\end{itemize}
\end{theorem}
In particular, the largest entries in the flag $h$-vector for Gorenstein* posets
are~$h_S$, for $S=\{0,2,4,\ldots\}$ and $S = \{1,3,5,\ldots\}$.
Readdy~\cite{Readdy-ext} showed that in the case of the crosspolytope, the 
maximum $h_S$ occurs only for these sets.

For the specific case of the simplex (Boolean lattice), 
Mahajan~\cite{Mahajan-Boolean} looked at 
inequalities among the coefficients of the $cd$-index.  He found, for 
example, that for the simplex
the coefficient of any $cd$-word of the form  $udv$ is greater
than or equal to the coefficient of $uccv$.  Furthermore the maximum
coefficient is, for $n$ even, the coefficient of $cd^jc$ 
with $j=(n-2)/2$ and, for
$n$ odd, the coefficient of $cdcd^jc$ (and that of
$cd^jcdc$, which is the same) with $j= (n-5)/2$.

Ehrenborg~\cite{lifting} gave a method for lifting any $cd$-inequality
to give inequalities in higher ranks.
For flag vectors of rational polytopes, a main source of inequalities was 
the nonnegativity of the $g$-vector (as in Definition~\ref{torich}) and
a form of lifting of these by convolution (Kalai~\cite{Kalai}).
Stenson~\cite{Stenson}
showed that the inequalities described in this section for $cd$-indices 
give flag vector inequalities that are not implied by the $g$-vector
convolution inequalities.

Murai and Yanagawa~\cite{Murai-Yana} defined {\em squarefree $P$-modules},
a generalization
of the Stanley--Reisner ring, and used it to generalize the $cd$-index to
a class of posets they call {\em quasi CW posets}.  They were then able to
prove that the coefficient of $w$ for a Gorenstein* poset is less than or
equal to the product of coefficients of associated $cd$-words having a
single $d$.  As a consequence, they get 
upper bounds on the $cd$-index of Gorenstein* posets when the number
$f_i$ of rank $i$ elements is fixed for all~$i$.

\section{Computing the $cd$-Index}
\subsection{Specific Polytopes and Posets}

Certain polytopes and posets have particularly nice $cd$-indices, often
connected to other combinatorial objects.
We will generally not define the associated combinatorial objects; the
interested reader can find details in the references.

Purtill's early results on nonnegativity of the $cd$-index 
(Theorem~\ref{Purtill}) resulted from studying CL-shellings of polytopes
\cite{Purtill}.
In this study he showed that the $cd$-index of the simplex is the
(noncommutative) Andr\'{e} polynomial of 
Foata and Sch\"{u}tzenberger~\cite{Foata-Schutz};
the Andr\'{e} polynomial is a generating function
for permutations satisfying certain descent properties. 
Purtill also extended the notion of Andr\'{e} permutations to signed 
permutations,  defined signed Andr\'{e} polynomials, and showed that
the signed Andr\'{e} polynomial is the $cd$-index of the crosspolytope
(the dual of the cube).  Note that reversing the monomials in the
signed Andr\'{e} polynomial gives the $cd$-index of the cube.
Thus, in the case of the simplex, crosspolytope and cube, each 
coefficient in the $cd$-index can be computed by counting (signed)
permutations with certain descent patterns.

Subsequently, Simion and Sundaram~\cite{Sundaram} defined the {\em simsun}
permutations, also counted by the Andr\'{e} polynomials.
Hetyei~\cite{Hetyei-var} gave an alternative set of permutations, which
he called {\em augmented Andr\'{e} permutations}, that give the $cd$
coefficients for the cube (and thus for the crosspolytope).
Billera, Ehrenborg and Readdy~\cite{c2d} gave formulas
for the $cd$-indices of the simplex, cube, and crosspolytope, with 
summations over all permutations.
Ehrenborg and Readdy~\cite{EhrRea-major} applied the connection in
the opposite direction, and
used the $ab$-index of the simplex and crosspolytope to study the
major index of permutations and signed permutations.

Inequalities for the $cd$-index of zonotopes were given in 
Section~\ref{Monotonicity}.  
Billera, Ehrenborg and Readdy~\cite{BER-zon} showed that $n$-dimensional
zonotopes span the generalized Dehn--Sommerville space and that they 
generate as an Abelian group all integral polynomials of degree $n$ in
$c$ and $2d$.
Bayer and Sturmfels~\cite{Bayer-Sturm} showed that the 
flag vector of an oriented matroid is determined by the underlying
matroid.  Billera, Ehrenborg and Readdy~\cite{c2d} gave an
explicit formula for the $cd$-index of an oriented matroid in terms of
the $ab$-index of its lattice of flats.  In particular, this gives
formulas for the $cd$-indices of zonotopes and of essential hyperplane
arrangements.
Ehrenborg, Readdy and Slone~\cite{Ehr-Rea-Slo} extended this to 
affine and toric hyperplanes.
In another direction it  was extended to ``oriented interval greedoids'' by
Saliola and Thomas~\cite{Saliola}.

Ehrenborg and Readdy gave recursive formulas for the $cd$-index of the simplex
and the cube \cite{EhrRea-coprod}.
They also gave recursive formulas for the $cd$-index of the lattice of regions
of the braid arrangements $A_n$ and $B_n$ \cite{EhrRea-Dowling}.
Joji\'{c}~\cite{Jojic-weighted} then gave the $cd$-index of 
the lattice of regions of the arrangements $D_n$ in terms of those
of $A_n$ and $B_n$.

Hsiao~\cite{Hsiao} gave a general construction of a class of Gorenstein* posets,
based on a signed version of the construction of a distributive lattice from
the order ideals of a general poset.  For the resulting ``signed Birkhoff
posets'' he gave a combinatorial description of the $cd$-index in terms of
peak sets of linear extensions of the underlying poset.

Two other combinatorial computations of the $cd$-index of a simplex (Boolean
lattice) were given by 
Fan and He~\cite{FanHe-Boolean} (based on 
 methods from the Bruhat order (Section~\ref{Bruhat})), and by
Karu~\cite{KaruM} (counting certain integer-valued functions).

\subsection{Operations on Posets}\label{operations}
Among the tools used in the study of $cd$-indices are results about the
effect on the $cd$-index of various operations on posets.
The methods used to develop many of these involve the coproduct, introduced
by Ehrenborg and Readdy~\cite{EhrRea-coprod}; we postpone discussion of 
that until Section~\ref{algebras}.

The most straightforward effect on the $cd$-index occurs for the join of
two posets.
\begin{definition}{\em
Given graded posets $P$ and $Q$, the {\em
join} $P\ast Q$ of $P$ and $Q$ is the poset on the set 
$(P\setminus \{\hat{1}\})\cup(Q\setminus\{\hat{0}\})$ with $x\preceq y$
in $P\ast Q$ in the following cases: 
\begin{itemize}
\item $x\preceq y$ in 
$P\setminus \{\hat{1}\}$ 
\item $x\preceq y$ in  $Q\setminus\{\hat{0}\}$
\item $x\in P\setminus \{\hat{1}\}$ and $y\in Q\setminus\{\hat{0}\}$.
\end{itemize}
}\end{definition}
\begin{theorem}[\cite{Stanley-flag}] 
If $P$ and $Q$ are Eulerian posets, then so is $P\ast Q$, and
\begin{equation*}\Phi_{P\ast Q}(c,d) = \Phi_P(c,d)\Phi_Q(c,d).\end{equation*}
\end{theorem}

The {\em pyramid} of a poset $P$ is the Cartesian product 
$\Pyr(P)=P\times B_1$,
where $B_1$ is the two-element chain.
The {\em prism} of a poset $P$ is the ``diamond product,'' 
$\Prism(P) = P\diamond B_2 = 
(P\setminus\{\hat{0}\})\times (B_2\setminus\{\hat{0}\})\cup \{\hat{0}\}$,
where $B_2$ is the Boolean lattice on two elements.
The dual operation to the prism operation takes $P$ to $\Bipyr(P)$.
(The terms come from the polytope context.)
These operations produce Eulerian posets from Eulerian posets.
Ehrenborg and Readdy computed the effect of these operations on the
$cd$-index.
They expressed this in terms of a couple of derivations on $cd$-words.  We
show one set of formulas; for others see~\cite{EhrRea-coprod}.
Define a derivation $D$ on $cd$-words by $D(c)=2d$ and $D(d) = cd +dc$.
\begin{theorem}[\cite{EhrRea-coprod}]\label{cd-operations}
Let $P$ be an Eulerian poset.  Then
\begin{itemize}
\item $\displaystyle \Phi(\Pyr(P)) = \frac{1}{2}
      \left[ \Phi(P)\,c+ c\,\Phi(P)+D(\Phi(P)) \right]$
\item $\displaystyle \Phi(\Prism(P)) = 
       \Phi(P)\,c+ D(\Phi(P))$
\item $\displaystyle \Phi(\Bipyr(P)) = 
       c\,\Phi(P)+  D(\Phi(P))$
\end{itemize}
\end{theorem}
Ehrenborg and Readdy~\cite{EhrRea-coprod}
also described the effect on the $cd$-index
of other operations on polytopes: truncation at a vertex, gluing polytopes
together along a common facet (in particular, performing a stellar subdivision
of a facet), and taking a Minkowski sum with a segment. 
They also gave a formula for the $ab$-index (flag $h$-vector) of the 
Cartesian product of arbitrary polytopes.
Ehrenborg and Fox~\cite{Ehr-Fox} gave recurrences for the $cd$-index 
of the Cartesian product and free join of polytopes.
Slone~\cite{Slone} gave a lattice path interpretation of the diamond product
of two $cd$-words. 
Fox~\cite{Fox} extended this interpretation to  all $cd$-words.
Ehrenborg, Johnson, Rajagopalan and Readdy~\cite{EJRR} gave formulas
for the $cd$-index of the polytope resulting from cutting off a face
of a polytope and for the $cd$-index of the regular CW complex 
resulting from contracting a face of the polytope.
S.~Kim~\cite{Kim-S} showed how the $cd$-index of a polytope can be expressed
when a polytope is split by a hyperplane.
Wells~\cite{Wells} generalized the idea of bistellar flips to (polytopal)
PL-spheres and computed the effect on the $cd$-index.

For a fixed graded poset $P$, one can form the poset $I(P)$ of all closed
intervals ordered by inclusion. Joji\'{c}~\cite{Jojic-int}
studied this poset, showed that
if $P$ is Eulerian then $I(P)$ is Eulerian, and computed the $cd$-index
of $I(P)$ in terms of that of $P$.
Joji\'{c} also computed the effect on the $cd$-index of the 
``$E_t$-construction'' of Paffenholz and Ziegler~\cite{Paffenholz}.

Hetyei~\cite{Hetyei-Tch} introduced the Tchebyshev transform on posets.
He used it to 
construct a sequence of Eulerian posets (one in each
rank) with a very simple formula for the $cd$-coefficients.  The $ce$-index
(a variation of the $cd$-index) of this poset is equivalent to the 
Tchebyshev (Chebyshev) polynomial.
Ehrenborg and Readdy~\cite{EhrRea-Tcheb} continued the study of the
Tchebyshev transform on general graded posets.
They showed that the Tchebyshev transform (of the first kind) preserves these
poset properties: Eulerian, EL-shellable and Gorenstein*.
In the Eulerian case, they computed the $cd$-index of 
$T(P)$ in terms of that of $P$ and showed that nonnegativity of the 
$cd$-index of $P$ implies nonnegativity of the $cd$-index of $T(P)$.
They showed that a second kind of Tchebyshev transform is a Hopf algebra
endomorphism on the Hopf algebra of quasisymmetric functions (see
Section~\ref{algebras}).

\subsection{Shelling Components}
Stanley~\cite{Stanley-flag} decomposed the $cd$-index of an \mbox{$n$-dimensional}
simplex (Boolean lattice) into parts based on a shelling of the simplex,
and used the parts for a formula for the $cd$-index of a simplicial 
Eulerian poset.
A simplicial Eulerian poset is an Eulerian poset such that for every 
$x\ne \hat{1}$, the interval $[\hat{0},x]$ is a Boolean lattice.
Note that the $h$-vector of a simplicial Eulerian poset is defined by
the transformation from the $f$-vector (mentioned in Section~\ref{Basic})
for simplicial polytopes.

Let $\sigma_0$, $\sigma_1$, \ldots, $\sigma_n$ be any ordering of the
facets of the $n$-simplex $\Sigma^n$; it is a shelling order.  Let 
$\hat{\Phi}_i^n(c,d)$ be the contribution to the $cd$ index of $\Sigma^n$ from
the faces added when $\sigma_i$ is shelled on.  
(For details see \cite{Stanley-flag}.)  We refer to these as the
{\em shelling components} of the $cd$-index.
\begin{theorem}

\mbox{}

\begin{itemize}
\item For all $i$, $0\le i\le n-1$, $\hat{\Phi}_i^n(c,d)\ge 0$
\item If $P$ is a simplicial Eulerian poset of rank $n+1$ with $h$-vector
      $(h_0,h_1,\ldots, h_n)$, then
      $\Phi_P(c,d) = \sum_{i=0}^{n-1} h_i\hat{\Phi}_i^n(c,d)$.
\end{itemize}
\end{theorem}  
As a consequence, Stanley proved the nonnegativity of the $cd$-index
for Gorenstein* simplicial posets before Karu's proof for general
Gorenstein* posets.
Stanley conjectured, and Hetyei~\cite{Hetyei-var} proved formulas for
the shelling components $\hat{\Phi}_i^n(c,d)$ in terms of Andr\'{e}
permutations and in terms of simsun permutations.
Ehrenborg and Readdy~\cite{EhrRea-coprod} gave a compact recursion
for these shelling components, and
Ehrenborg~\cite{ehr-shell} gave more recursions for them.

Ehrenborg and Hetyei~\cite{EhrHet-cub}
developed the analogous results for cubical Eulerian
posets, that is, Eulerian posets whose lower intervals are isomorphic to
the face lattice of a cube.
Billera and Ehrenborg~\cite{Bil-Ehr}
gave a formula for the contribution of each facet in a shelling of a polytope.
Lee~\cite{Lee-sweeping} described a dual approach: the calculation of the 
$cd$-index by ``sweeping'' a hyperplane through the polytope, keeping track
of the contribution at each vertex.  

\section{Bruhat Order}\label{Bruhat}
The original motivation for the $cd$-index came from the combinatorial study
of convex
polytopes, but Reading began the study of the $cd$-index for another important
class of Eulerian posets: intervals in the Bruhat order of Coxeter groups.  
In short, for $v$ and $w$ elements of a Coxeter group, $v\prec w$ if and only
if some reduced word representation of $v$ is a subword of a reduced
word for $w$.
An interval in the Bruhat order is Eulerian and shellable, and hence
Gorenstein*.
For example, the $cd$-index of the Bruhat order of ${\mathcal S}_4$ is
$c^5+c^3d+2c^2dc+2cdc^2+dc^3+2cd^2+dcd+2d^2c$.
For more information on Bruhat order in our context, see 
\cite{Bjorner-Wachs,Reading}.

Reading~\cite{Reading} gave a recursive formula for the $cd$-index of a 
Bruhat interval.  He showed that Bruhat intervals span the 
generalized Dehn--Sommerville space, and gave an explicit basis.
\begin{theorem}
The set of $cd$-indices of Bruhat intervals of rank $r$ spans the
affine span of $cd$-indices of Eulerian posets of rank $r$.
\end{theorem}

The Bruhat order of a universal Coxeter group contains intervals isomorphic
to the face lattices of certain polytopes, the duals of stacked polytopes.
Reading conjectured that these have maximum $cd$-indices.
\begin{conjecture}[\cite{Reading}]
Let $(W,S)$ be a Coxeter system, and let  $[u,v]$ be an interval in the Bruhat
order of $W$ with $u$ of length $k$ and $v$ of length $n+k+1$.
Then the $cd$-index of $[u,v]$ is coefficientwise less than or equal to
the $cd$-index of a dual stacked $n$-polytope with $n+k+1$ facets.

In particular the $cd$-index of an interval $[1,v]$ is less than or
equal to the $cd$-index of a Boolean lattice.
\end{conjecture}

The $cd$-index can be found in the peak algebra of quasisymmetric functions
\cite{BHvW} (see Section~\ref{algebras}).
Billera and Brenti~\cite{BilBrenti} used this to define for Bruhat intervals
the {\em complete $cd$-index}, a nonhomogeneous polynomial in $c$ and $d$,
whose homogeneous part of top degree is the $cd$-index.
They used this to give an explicit computation of the Kazhdan--Lusztig
polynomials of the Bruhat intervals for any Coxeter group.
They conjectured that all coefficients of the complete $cd$-index are
nonnegative for all Bruhat intervals.

Besides the top degree terms, whose nonnegativity follows from Karu's 
theorem, the nonnegativity of certain coefficients in the complete $cd$-index
of Bruhat intervals have been verified 
\cite{Blanco-shortest,FanHe-nonneg,Karu-Bruhat}.
Blanco~\cite{Blanco-complete}
used CL-labeling due to Bj\"{o}rner and Wachs~\cite{Bjorner-Wachs} to
describe the computation of the complete $cd$-indices of dihedral Bruhat
intervals (those isomorphic to intervals in a dihedral reflection
subgroup) and Bruhat intervals in universal Coxeter groups.

Blanco~\cite{Blanco-shortest} defined the shortest path poset in a
Bruhat interval, and showed that if the poset has a unique maximal
rising chain then it is a Gorenstein* poset.
Y.~Kim~\cite{Kim-Y} studied the uncrossing partial order of matchings on $[2n]$,
which is isomorphic to a subposet of the dual Bruhat order of 
affine permutations.  
He gave a recursion for the $cd$-indices of intervals in this poset.

\section{Algebras}\label{algebras}
The $h$-vector of a simplicial polytope and the flag $h$-vector of Eulerian
posets have interpretations in the Stanley--Reisner ring of the polytope or
of the order complex of the poset.
The $cd$-index is not found naturally in this ring.  It turns out that other
algebras are better habitats for the $cd$-index.
There is an extensive literature on these algebras, and this survey will
only touch the surface.  For a deeper look, the reader is directed to
(in chronological order)
\cite{Ehr-Hopf,EhrRea-coprod,Billera-Liu,Bergeron,Aguiar,Newtonian,BHvW,
BilBrenti,Grujic,KaruM,Brenti-Cas}.
See Billera~\cite{Bil-ICM} for a survey of some of these connections.

Perhaps the beginning of the story is \cite{Ehr-Hopf}, where Ehrenborg
gave a Hopf algebra homomorphism from the Hopf algebra of posets (the ``reduced
incidence Hopf algebra'') to the Hopf algebra of quasisymmetric functions.
The homomorphism gives some (not all) known results on flag vectors of 
Eulerian posets.

An important concept underlying some of the algebraic structures is
the convolution of flag numbers, introduced by Kalai~\cite{Kalai}.
The entries of the flag vector, $f_S^n$ are considered as 
functions from rank $n$ graded posets to nonnegative integers.  
A convolution product is defined:  for $P$ a poset of rank $n+m$,
\begin{equation*}f_S^n\ast f_T^m(P)= \sum_{\genfrac{}{}{0pt}{}{x\in P}{\rho(x) = n}}
f_S^n([\hat{0},x])f_T^m([x,\hat{1}])= 
f_{S\cup\{n\}\cup(T+n)}^{n+m}(P).\end{equation*}
Ehrenborg and Readdy~\cite{EhrRea-coprod} described a coproduct on the
vector space spanned by graded posets by
\begin{equation*}\Delta(P) = \sum_{\genfrac{}{}{0pt}{}{x\in P}{\hat{0}<x<\hat{1}}} [\hat{0},x]\otimes
[x,\hat{1}],\end{equation*} and a coproduct on the noncommutative polynomial ring
$k\langle a,b\rangle$ by 
\begin{equation*}\Delta(u_1u_2\cdots u_n) = \sum_{i=1}^n u_1\cdots u_{i-1}\otimes
u_{i+1}\cdots u_n .\end{equation*}
(Recall that the $ab$-polynomial of a graded poset has as its coefficients
the flag \mbox{$h$-numbers}, and these coefficients can be written in terms of the
flag $f$-vector.)
They showed that the $ab$-index is a Newtonian coalgebra map between the 
resulting coalgebras,
and showed that this map takes the subalgebra spanned by Eulerian posets to
the subalgebra $k\langle c,d\rangle$.
They used this to derive the formulas for the effect of various operations 
on the $cd$-index (see Section~\ref{operations}).

Billera and Liu~\cite{Billera-Liu} introduced a graded algebra of flag
operators on posets.  (In \cite{BHvW} it is referred to as the
``algebra of forms on Eulerian posets.'')  The flag number functions
$f_S^n$  span a graded vector space over $\bf Q$, $A = \oplus_{n\ge 0}A_n$,
with \mbox{$A_n = \{ \sum_{S\subseteq [n-1]} \alpha_S f_S^n: \alpha_S\in {\bf Q}\}$}.
With the convolution product, $A$ becomes a graded algebra, and can be generated
by the trivial flag operators $f_\emptyset^j$.
Billera and Liu determined the two-sided ideal of $A$ of elements that vanish 
for all Eulerian posets.  
Write $A_{\mathcal E}$ for  the quotient of $A$ by this ideal.
\begin{theorem}[\cite{Billera-Liu}]

\mbox{}

\begin{itemize}
\item
As graded algebras, $A\cong {\bf Q}\langle y_1,y_2,\ldots\rangle$, the free
graded associative algebra on generators $y_i$ of degree $i$.
The isomorphism is determined by $f_\emptyset^j \mapsto y_j$.
\item As graded algebras, 
$A_{\mathcal E}\cong {\bf Q}\langle y_1,y_3, y_5,\ldots\rangle$.
\end{itemize}
\end{theorem}

The $ab$-polynomial of an Eulerian poset is in the 
(noncommutative) polynomial ring $A_{\mathcal E}\langle a,b\rangle$.  
Billera and Liu gave another proof of the existence of the $cd$-index
for Eulerian posets by proving
\begin{theorem}[\cite{Billera-Liu}]
As a polynomial with coefficients in $A_{\mathcal E}$, the $ab$-polynomial of
every Eulerian poset is in $A_{\mathcal E}\langle c,d\rangle$.
\end{theorem}
The authors also described one-sided ideals of the graded algebra $A$ 
representing flag vector relations on simplicial polytopes and on cubical 
polytopes, and gave dimension arguments from the resulting $A$-modules.

Next we look at the {\em peak algebra}, introduced by 
Stembridge~\cite{Stembridge} in the study of enriched $P$-partitions, and described
in terms of the $cd$-index in \cite{BHvW}, following~\cite{Bergeron}.
For $T=\{t_1,t_2,\ldots, t_k\}\subseteq [n]$, let 
\begin{equation*}M_T = \sum_{1\le i_1<i_2<\cdots <i_k} 
x_{i_1}^{t_1}x_{i_2}^{t_2-t_1}x_{i_3}^{t_3-t_2}\cdots x_{i_k}^{t_k-t_{k-1}}.
\end{equation*}
This is a quasisymmetric function in $x_i$ ($i\ge 1$); the sum is over all 
increasing $k$-tuples.
For a $cd$-word $w = c^{n_1}dc^{n_2}d\cdots c^{n_k}dc^m$ of degree $n$,
define $S_w\subseteq [n]$ by
$S_w = \{n_1+2,n_1+n_2+4,\ldots, n_1+n_2+\cdots + n_k+2k\}$
\begin{definition}[from \cite{BHvW}]{\em
Let $\mathcal Q$ be the algebra of quasisymmetric
functions over $\bf Q$ in the variables $x_1, x_2, \ldots$.
The {\em peak algebra} is the subalgebra $\Pi$ of $\mathcal Q$ generated by the 
elements 
$\displaystyle \Theta_w=\sum_{\genfrac{}{}{0pt}{}{T}{S_w\subseteq T\cup (T+1)}} 2^{|T|+1}M_T$
for each $cd$-word $w$.
}\end{definition}

The two algebras have natural coproducts that make them Hopf algebras.
\begin{theorem}[\cite{Bergeron}]
The algebra of forms on Eulerian posets and the peak algebra are dual Hopf
algebras.
\end{theorem}

Following \cite{Ehr-Hopf}, Billera, Hsiao and van Willigenburg~\cite{BHvW} 
considered the quasisymmetric representation of the flag $f$-vector, 
$F(P) = \sum_{S\subseteq[n]}f_S(P)M_S$, and showed that when written in terms
of Stembridge's basis $\{\Theta_w\}$ the coefficients are the $cd$-index of
$P$ (modified by factors of 2).
Then an Eulerian poset has a nonnegative $cd$-index when its
peak quasisymmetric function has a nonnegative representation in terms
of this basis.

Aguiar~\cite{Aguiar} took an algebraic approach to constructing the $ab$-index
of a
poset using a morphism from the algebra of all posets to the noncommutative
polynomial algebra $k\langle a,b\rangle$, considered as infinitesimal
Hopf algebras.  This perspective enabled him to consider generalizations
of the $ab$-index to weighted posets and the relative $ab$-index.
Using a map involving the zeta and M\"{o}bius functions of a poset,
he found the infinitesimal Hopf algebra corresponding to Eulerian posets, 
that is, the $cd$-index.

Ehrenborg and Readdy~\cite{Newtonian} gave yet another proof of the existence
of the $cd$-index for Eulerian posets through the homology of their
Newtonian coalgebra generated by $a$ and $b$, mentioned above.

Karu~\cite{KaruM} gave an algebraic formulation of a conjecture of
Murai and Nevo~\cite{MuraiNevo-cd} and proved it for a special case.
The statement depends on a representation of $cd$-words by $0/1$-vectors:
$\mbox{mdeg}(w)$ is the $0/1$-vector obtained from $w$ by replacing each
$c$ by 0 and each $d$ by 10.  For $w$ a $cd$-word and $v=\mbox{mdeg}(w)$,
write $\Phi_{P,v}$ for the coefficient of $w$ in $\Phi_P(c,d)$.
For vectors $v\in {\bf Z}^n$ that are not equal to $\mbox{mdeg}(w)$ 
for all $cd$-words $w$, let $\Phi_{P,v}=0$.
\begin{theorem}[\cite{KaruM}]
Let $P$ be the poset of a Gorenstein* simplicial complex of dimension $n$.
Then there exists a standard ${\bf Z}^n$-graded $k$-algebra 
$A=\oplus_v A_v$ such that $\Phi_{P,v}= \dim A_v$.
\end{theorem}
Murai's and Nevo's conjecture was that this theorem holds for all 
Gorenstein* posets.

Fine \cite{Fine-arxiv} suggests a successor to the $cd$-index. 
Like the $cd$-index, it provides a special basis for the flag vector ring,
which he calls $\mathcal{R}$. The successor is a basis that, 
conjecturally, has this property: the structure coefficients for product and
pyramid (which he calls cone) are all nonnegative integers.
See also Theorem~\ref{cd-operations}
and~\cite{Ehr-Fox} for product and pyramid in the $cd$-index.
If $\mathcal{R}$ were the representation ring of some algebraic object
$\mathcal{G}$ that satisfies Schur's lemma, then the existence of such a basis
would be immediate.

\section{Related Parameters}
Here we mention some generalizations of the $cd$-index to broader settings.

Ehrenborg  and Readdy~\cite{EhrRea-cubical} considered a generalization of
the lattice of the cube. 
For various $r$, take the 
posets $M_r$ with $r$ minimal elements and a unique maximum element,
then take the Cartesian product and add a unique minimum element.
The result they called an $\bf r$-cubical lattice, where $\bf r$ is
the vector of the various values of $r$ in the factors.  They defined
a generalized $cd$-index for these lattices, and showed that the
coefficients count a class of permutations generalizing the Andr\'{e}
permutations.

Ehrenborg~\cite{kEulerian}  considered a relaxation of the Eulerian condition.  
A poset $P$ is {\em $k$-Eulerian} if every interval of rank $k$ is
Eulerian.  He showed that the $ab$-index of a $k$-Eulerian poset can
be written in terms of $c$, $d$ and $e^{2k+1} = (a-b)^{2k+1}$.
He also related the $k$-Eulerian posets to an ideal in the Newtonian
coalgebra.

Ehrenborg, Hetyei and Readdy~\cite{EHR-level} considered {\em Eulerian level
posets}, infinite posets with a certain uniformity at each level and whose
finite intervals are Eulerian.  For these posets they extended the 
$cd$-polynomial to a $cd$-series.

As mentioned in Section~\ref{other-ineq}, Murai and Yanagawa~\cite{Murai-Yana}
considered a class of ``quasi CW posets'' and defined an extended $cd$-index 
for this class.
Gruji\'{c} and Stojadinovi\'{c}~\cite{Grujic} developed an analogue of
the $cd$-index for ``building sets'' (see~\cite{DeConcini}) by
introducing a Hopf algebra of building sets and 
mimicking the known Hopf algebra construction of the $cd$-index.

Ehrenborg, Goresky and Readdy~\cite{EGR}
extended the definition of the $cd$-index to ``quasi-graded posets''
with a generalization of the notion of Eulerian.  These posets arise
from Whitney stratified manifolds.  
They studied in particular the $cd$-indices
of semisuspensions.
Ehrenborg and Readdy~\cite{EhrRea-Manifold} applied this to study
manifold arrangements.

Dornian, Katz and Tsang~\cite{DKT} introduced the ``mixed $cd$-index''
of a strong formal subdivision of an Eulerian poset.  This is defined in terms
of the local $cd$-index of Karu \cite{Karu-cd}, and maps to the 
``mixed $h$-polynomial'' of \cite{KS}, a multivariate polynomial invariant of 
subdivisions of polytopes.

Murai and Nevo~\cite{MuraiNevo-cd} related the $cd$ index of a class of
Eulerian posets to the $f$-vector of a simplicial complex.  
The relation is via the $\gamma$-vector, first introduced by Gal~\cite{Gal}.
Here the set of $S^*$-shellable spheres is a subset of Stanley's $S$-shellable
spheres, and includes convex polytopes.
\begin{theorem}[\cite{MuraiNevo-cd}]
Let $P$ be an $(n-1)$-dimensional $S^*$-shellable regular CW sphere, with 
$cd$-index $\Phi_P(c,d)$.  Define $\delta_i$ by 
$\Phi_P(1,d) = 
\delta_0+\delta_1d+\cdots + \delta_{\lfloor n/2\rfloor}d^{\lfloor n/2\rfloor}$.
Then there exists an $\lfloor n/2\rfloor$-colored simplicial complex
$\Delta$ such that $\delta_i = f_{i-1}(\Delta)$ for 
$0\le i\le\lfloor n/2\rfloor $.
\end{theorem}

This is a survey of the $cd$-index, but as mentioned in Section~\ref{Basic},
of great interest in the study of Eulerian posets and, in particular, of 
convex polytopes, is the toric $h$-vector.  The toric $h$-vector contains
much less information than the $cd$-index.  However, Lee was led by his 
study of sweeping a hyperplane through a polytope to
an extension of the toric $h$-vector that is equivalent to the $cd$-index
\cite{Lee-sweeping}.

\section{Conclusion}
The introduction of the $cd$-index opened up many directions of research on
Eulerian posets.
There are many specific open questions, but the overriding
issue is to find combinatorial interpretations for the coefficients, beyond
those in special cases mentioned here.
I hope that this survey will become outdated soon, because of  
significant research advances.

\section*{Acknowledgments}

The referees were generous with corrections and suggestions, which improved
the paper.  I also wish to thank Louis Billera, Jonathan Fine, Eran Nevo
and Richard Stanley, who gave helpful feedback on a draft of this paper.

\providecommand{\bysame}{\leavevmode\hbox to3em{\hrulefill}\thinspace}
\providecommand{\MR}{\relax\ifhmode\unskip\space\fi MR }
\providecommand{\MRhref}[2]{%
  \href{http://www.ams.org/mathscinet-getitem?mr=#1}{#2}
}
\providecommand{\href}[2]{#2}

\end{document}